\documentclass[12pt]{article}
\begin{document}
\newcommand{\qed}{\ \
\mbox{\rule{8pt}{8pt}}\vspace{0.3cm}\newline}
\newcommand{\ia}{{\bf I}_{A_i}}
\newcommand{\ba}{\widetilde{{\bf KC}}_i}
\newcommand{\bba}{{\bf KC}_i}
\newcommand{\ra}{\longrightarrow}
\newcommand{\pe}{{\cal P}}
\newcommand{\dep}{{\cal DP}}
\newcommand{\adp}{{\cal DP}_d^{af}}
\newcommand{\ot}{\otimes}
\newcommand{\ser}{{\bf S}}		
\newcommand{\la}{\lambda}
\newcommand{\rec}{\raisebox{-1ex}{\ $\stackrel{\textstyle{\stackrel{\textstyle{\longleftarrow}}{\longrightarrow}}}{\longleftarrow}$\ }}
\title{On Serre functor in the category of strict polynomial functors}
\author{Marcin Cha\l upnik
\thanks{The author was supported by the  grant (NCN) 2011/01/B/ST1/06184.}\\
\normalsize{Institute of Mathematics, University of Warsaw,}\\
\normalsize{ul.~Banacha 2, 02--097 Warsaw, Poland;}\\
\normalsize{e--mail: {\tt mchal@mimuw.edu.pl}}}
\date{\mbox{}}
\newtheorem{prop}{Proposition}[section]
\newtheorem{cor}[prop]{Corollary}
\newtheorem{theo}[prop]{Theorem}
\newtheorem{lem}[prop]{Lemma}
\newtheorem{defi}[prop]{Definition}
\newtheorem{defipro}[prop]{Definition/Proposition}
\newtheorem{fact}[prop]{Fact}
\newtheorem{exa}[prop]{Example}
\newcommand{\ka}{{\mbox {\bf k}}}
\newcommand{\ca}{\mbox{$\cal{A}$}}
\maketitle
\begin{abstract}
We introduce and study a Serre functor in the category ${\cal P}_d$ of strict polynomial functors over a field of positive characteristic. By using it we obtain  
the Poincar\'e duality formula for Ext--groups from [C3] in  elementary way. We also show that the derived category of the category of affine strct polynomial functors in some cases carries the structure of Calabi-Yau category. 
\\\mbox{}\vspace{0.1cm}
\end{abstract}
\section{Introduction}
In the present article we study a Serre functor in the category ${\cal P}_d$ of strict polynomial functors of degree $d$ over a field of positive characteristic.  Although the existence of a Schur functor
in our context follows from general theory, its interplay with various  structures living on ${\cal P}_d$ (Frobenius twist, affine subcategories, blocks)
has some interesting consequences.\newline
Our paper can be naturally divided into two parts: Section 2 and Sections 3--5.
\newline The objective of the first part is to give an elementary and self--contained account of the Serre functor $\ser$ in ${\cal P}_d$
and to quickly obtain with the aid of $\ser$  the Poincar\'e duality formula   from [C3]:\vspace{0.5cm}\newline
{\bf Corollary 2.4}	 Let $\la$ be a  Young diagram of weight $d$ which is single in its block (we call such a diagram and its block {\em basic}),
let $\mu$ be any Young diagram of weight
	$p^id$.  Let $F_{\la}, F_{\mu}$ be the corresponding simple objects. Then
	\[ \mbox{Ext}^s_{{\cal P}_{dp^i}}(F_{\la}^{(i)}, F_{\mu})\simeq \mbox{Ext}^{2d(p^i-1)-s}_{{\cal P}_{p^id}}(F_{\la}^{(i)}, F_{\mu})^*.\]
We belive that the approach to the Poincar\'e duality presented here is more
intuitive than that taken in [C3] and  that it is also better adapted for possible generalizations. \newline
Then in Sections 3--5 we turn  attention to the category ${\cal P}_d^{af_i}$ of $i$--affine strict polynomial functors of degree $d$
and we introduce a (somewhat weaker version of) Serre functor on its derived
category  ${\cal DP}_d^{af_i}$. The main goal of this part of the paper
is to put the formula from Corollary 2.4 into a wider categorical context. 
In general, a Serre functor produces  Poincar\'e duality in Ext--groups when it acts on some
object as the shift functor. Indeed, we see in our Proposition 2.3 that this exactly
happens for some Frobenius twisted strict polynomial functors (in fact, Corollary
2.4 is a formal conequence of Proosition 2.3). We provide a
categorical interpretation of this phenomenon by finding certain subcategories
of ${\cal DP}_{dp^i}$ on which the Serre functor is isomorphic to the shift functor  (such categories are
called Calabi--Yau). In fact, we define these Calabi--Yau categories as the images of the
basic blocks from ${\cal P}_d$ in ${\cal DP}_d^{af_i}$ (we call these subcategories ``basic (affine) semiblocks''). Appearing of
${\cal DP}_d^{af_i}$ here is quite natural, since it is a full triangulated subcategory
of ${\cal DP}_{dp^i}$ generated by the Frobenius twisted objects [C4, Theorem~5.1]. 
Thus we have succeeded in providing a categorical interpretation of the both
assumptions in Corollary 2.4: we specialize to ${\cal DP}_d^{af_i}$ because 
$F_{\la}^{(i)}$ is twisted and we restrict to the image of the block containing 
$F_{\la}$ to take advantage of the fact that $\la$ is basic. \newline Below we describe the contents of Sections 3--5 in some detail  since the considerations there are much
more involved than those in Section 2.\newline
 In Section 3 we study a Serre functor  in  ${\cal DP}_d^{af_i}$. We start
with reviewing basic properties of the categories ${\cal P}_d^{af_i}$ and ${\cal DP}_d^{af_i}$. This is mainly recollecting some facts from [C4]  where these concepts were introduced and adapting them to a slightly more general setting
of  ``multiple twists'' in which we work in the present article. Then we proceed
to define the Serre functor $\ser^{af_i}$ in ${\cal DP}_d^{af_i}$. However, we need to adapt this
notion to the fact that  ${\cal DP}_d^{af_i}$ has infinite dimensional Hom--spaces.
Hence, technically, we define ``a weak Serre functor '' (Definition 3.3) on
${\cal DP}_d^{af_i}$. \newline In Section 4 we introduce  ``the semiblock decomposition'' of   ${\cal DP}_d^{af_i}$. This is a collection of reflective subcategories of
${\cal DP}_d^{af_i}$ indexed by the set of blocks in ${\cal P}_d$. They generate
${\cal DP}_d^{af_i}$ but in contrast to genuine blocks  they are not orthogonal.
We believe that this structure deserves further investigation, in particular we conjecture that
the semiblocks form a set of strata of certain stratification of ${\cal DP}_d^{af_i}$. In the present article we restrict ourself to introducing the affine derived Kan 
extension and Serre functor on the semiblocks. This is a non--trivial
task due to the non--orthogonality of semiblocks.\newline
In Section 5 we focus on the basic semiblocks, i.e. the subcategories of
${\cal DP}_d^{af_i}$ which correspond to the blocks containing a single simple object. We show (Theorem 5.1)
that they are Calabi--Yau categories thus providing the promised categorical interpretation
of  our Poincar\'e duality formula. We finish our paper by giving various 
explicit descriptions of basic semiblocks as  categories of DG--modules over certain graded algebras (Proposition 5.5, Corollary 5.6) which should make them  easier to handle. \newline Hence
the main result of this part of the paper is Theorem 5.1 which shows that the basic semiblocks are
Calabi--Yau categories. This is also important from a more general point of view. 
Namely, it seems that the basic semiblocks constitute sort of building blocks for 
${\cal DP}_{dp^i}$, since many homological problems concerning strict 
polynomial functors can be reduced to statements about basic semiblocks.
 For example, the classical line of research started in
[FS] and [FFSS] may be thought of as expanding of our understanding of the basic
semiblock ${\cal DP}_1^{af_i}$ onto the whole   ${\cal DP}_{p^i}$. To be more specific, let us consider
the fundamental problem of computing the Ext--groups between simple
objects in ${\cal P}_{dp^i}$. 
Then one can hope that by using tools like the the (Schur)--de Rham complex this problem can be reduced to 
that of computing 
\[\mbox{Ext}^*_{{\cal P}_{dp^i}}(F_{\la}^{(i)},F_{\mu})\] 
 where $\la$ is a basic Young diagram. This computation can be transferred
 by  the affine derived Kan extension to the basic semiblock containing $\la$. Therefore a better understanding of the internal structure of the basic semiblocks seems to be
 an important step towards understanding ${\cal DP}_{dp^i}$ in general. 
\section{Serre functor in ${\cal DP}_d$}
Let ${\cal P}_d$ be the category of strict polynomial functors of degree $d$ over
a fixed field $\ka$ of characteristic $p>0$ as defined in [FS]. For a finite dimensional $\ka$--vector space $U$ we define the strict polynomial functor $S^d_{U^*}\in {\cal P}_d$ by the formula
\[V\mapsto S^d(\mbox{Hom}(U,V)).\]
We recall from [FS, Th.~2.10] the natural in $U$ isomorphism
\[\mbox{Hom}_{{\cal P}_d}(F,S^d_{U^*})\simeq F(U)^*\]
	for any $F\in{\cal P}_d$.
In fact, when we interpret ${\cal P}_d$ as a functor category as is done e.g. in 
[FP, Sect.~3], this formula is just the Yoneda lemma. Hence later on we will refer to
this formula as to the Yoneda lemma. It immediately follows from this formula
that $S^d_{U^*}$ is injective and it was shown in [FS, Th.~2.10] that
if $\dim(U)\geq d$ then $S^d_{U^*}$ is a cogenerator of ${\cal P}_d$.
Dually, we have a family of projective objects $\Gamma^d_{U^*}$ for which
the Yoneda lemma gives the isomorphism
\[\mbox{Hom}_{{\cal P}_d}(\Gamma^d_{U^*},F)\simeq F(U)\]
for any $F\in{\cal P}_d$.
\newline
Let ${\cal DP}_d$ denote the bounded derived category of ${\cal P}_d$.
We would like to define a Serre functor in the sense of [BK] on ${\cal DP}_d$.
The problem of existence of Serre functor on a triangulated category is well
understood [BV, RV]. 
The existence of Serre functor on ${\cal DP}_d$ follows from the fact that it is
equivalent to the bounded derived category of  category of finitely generated modules over a finite dimensional algebra of finite homological dimension. 
When we translate the construction from [BV] into the context of functor categories we get
\begin{defi}
We define a  functor $\ser: \dep_d\ra\dep_d$ by the formula
\[\ser(F)(V):=\mbox{Hom}_{\dep_{d}}(S^d_{V^*},F).\]
\end{defi}
This functor was also mentioned in [Kr]. In our article we start a systematic study of its properties. Our exposition is quite elementary, independent of generalities of [BV]
and self--contained (with exception of a few places where we refer to [C1] to
avoid repeting the same arguments). We start by warning the reader that we chose
to work with less common left Serre functors which are easier to describe in the
framework of functor categories (although also the right Serre functor can be explicitly 
defined by using either the Kuhn duals, as we do in the proof of Theorem 2.2.4, or  the monoidal structure on ${\cal P}_d$ introduced in [Kr]). We gather below
basic properties of $\ser$. Parts 1, 4, 5 essentially follow from generalities, parts
2 and 3 are more specific to ${\cal P}_d$. 
\begin{theo}
The functor $\ser$ satisfies the following properties
\begin{enumerate}
	\item There is a natural in $U$ isomorphism in ${\cal P}_d$
	\[\ser(S^d_{U^*})\simeq\Gamma^d_{U^*}.\]
	\item  There is an isomorphism of functors
	\[\ser\simeq \Theta\circ \Theta\]
	where $\Theta$ is the ``Koszul duality'' functor from [C1].
   \item For any $F\in \dep_d$, $G\in\dep_{d'}$ there are isomorphisms in respectively
   $\dep_{dp^i}$, $\dep_{d+d'}$
    	\begin{itemize}
    	\item $\ser(F^{(i)})\simeq\ser(F)^{(i)}[-2d(p^i-1)]$
        \item $\ser(F\ot G)\simeq \ser(F)\ot \ser(G).$
        \end{itemize}   
    \item \ser\ is a self-equivalence of $\dep_d$.
    \item There is a natural in $F,G\in\dep_d$ isomorphism 
    \[{\mbox Hom}_{\dep_d}(F,G)\simeq{\mbox Hom}_{\dep_d}(\ser(G),F)^*,\]
 that is, \ser\ is a left Serre functor in the sense of [BK].
\end{enumerate}
\end{theo}
{\bf Proof:}  To see the first part we recall that since $S^d_{U^*}$ is injective,
we have  
\[\ser(S^d_{U^*})(V)=\mbox{Hom}_{\dep_{d}}(S^d_{V^*},S^d_{U^*})\simeq
\mbox{Hom}_{\pe_{d}}(S^d_{V^*},S^d_{U^*})\simeq S^d_{V^*}(U^*)^*\simeq
\Gamma^d_{U^*}(V)\]
 by the Yoneda lemma.\newline
In fact, [C1, Fact~2.2] can be easily extended to the ``parameterized version'':
\[\Theta((S_{\la})_{U^*})\simeq (W_{\widetilde{\la}})_{U^*}.\]  
From this we obtain the isomorphisms $\Theta(S^d_{U^*})\simeq 
\Lambda^d_{U^*}$ and
$\Theta(\Lambda^d_{U^*})\simeq \Gamma^d_{U^*}$ which give the second part. \newline
The formulae from part 3 follow from the analogous facts holding for $\Theta$ [C1, Fact~2.6].\newline
It is immediate that the``right Serre functor'' $\ser_l:=(-)^{\#}\circ\ser\circ(-)^{\#}$ where $(-)^{\#}$ is the Kuhn duality
is the inverse of \ser\ (c.f. [C1, Def. 2.3, Cor.~2.4]),
which gives the fourth part.\newline In order to obtain the last part, it suffices to establish
a natural in $U$ isomorphism
\[\mbox{Hom}_{\dep_d}(F,S^d_{U^*})\simeq
\mbox{Hom}_{\dep_d}(\ser(S^d_{U^*}),F)^*.\]
By the first part and injectivity of $S^d_{U^*}$ and projectivity  of 
$\Gamma^d_{U^*}$
it reduces to
\[\mbox{Hom}_{\pe_d}(F,S^d_{U^*})\simeq\mbox{ Hom}_{\pe_d}(\Gamma^d_{U^*},F)^*,\] 
which follows from the Yoneda  lemma. \qed
The fact that \ser\ is a Serre functor can be used to obtain the Poincar\`e like formulae for the Ext--groups, provided that we are able to compute $\ser(F)$ in some interesting 
cases. We shall ilustrate this idea by re--obtainng the most important example
of the Poincar\`e duality formula for Ext--groups in ${\cal P}_d$ established in 
[C3].\newline Let $\la$ be the Young diagram of weight $d$ which is a $p$--core.
We recall that the blocks in ${\cal P}_d$ are indexed by the $p$--core Young diagrams
of weight $d-jp$  and that the block labeled by $\la$ contains only one 
simple object $F_{\la}$. We call such a Young diagram $\la$  and the corresponding block {\em basic}.

\begin{prop}
Let $\la$ be a basic Young diagram. Then 
	\[\ser(F_{\la}^{(i)})\simeq F_{\la}^{(i)}[-2d(p^i-1)].\]
	\end{prop} 
	{\bf Proof:} Since $F_{\la}$ is single in its block, we have isomorphisms $F_{\la}\simeq S_{\la}\simeq W_{\la}$.
	Therefore
	\[\Theta(F_{\la})\simeq\Theta(S_{\la})=W_{\widetilde{\la}}.\]
	Now, since also $F_{\widetilde{\la}}$ is single in its block,  we obtain
	\[\Theta(W_{\widetilde{\la}})\simeq\Theta(S_{\widetilde{\la}})=W_{\la}
	\simeq F_{\la}.\]
	Thus we see that $\ser(F_{\la})\simeq F_{\la}$ and  our formula 
	follows from Theorem 2.1.3.\qed
	The Poincar\`e duality formula [C3, Example~3.3] is  a formal consequence
	of Proposition 2.3.
	\begin{cor}
	Let $\la$ be a basic Young diagram, $\mu$ be any Young diagram of weight
	$p^id$, and $F_{\la}, F_{\mu}$ be the corresponding simple objects. Then
	\[ \mbox{Ext}^s_{{\cal P}_{dp^i}}(F_{\la}^{(i)}, F_{\mu})\simeq \mbox{Ext}^{2d(p^i-1)-s}_{{\cal P}_{p^id}}(F_{\la}^{(i)}, F_{\mu})^*.\]
 \end{cor}
 {\bf Proof:} By applying the Kuhn duality (and using the fact that simple objects are self--dual) and then the Serre functor  we obtain:
\[\mbox{Ext}^s_{{\cal P}_{p^id}}(F_{\la}^{(i)}, F_{\mu})\simeq
\mbox{Ext}^s_{{\cal P}_{p^id}}(F_{\mu},F_{\la}^{(i)})=
\mbox{Hom}_{{\cal DP}_{dp^i}}(F_{\mu},F_{\la}^{(i)}[s])\simeq\]
\[\mbox{Hom}_{{\cal DP}_{dp^i}}(\ser(F_{\la}^{(i)}[s]),F_{\mu})^*\simeq
\mbox{Hom}_{{\cal DP}_{dp^i}}(F_{\la}^{(i)}[s-2d(p^i-1)], F_{\mu})^*\simeq\]
\[\mbox{Ext}^{2d(p^i-1)-s}_{{\cal P}_{dp^i}}(F_{\la}^{(i)}, F_{\mu})^*.
\]
\qed
In the next part of the paper we will describe a categorical phenomenon which
is responsible for turning the Serre duality into the Poincar\`e duality when
one deals with the Frobenius twists of strict polynomial functors.  
\section{Serre functor for affine functors}
\subsection{Review of $i$--affine functors} 
The category of affine strict polynomial functors ${\cal P}_d^{af}$ was studied
in [C4]. In the present paper we introduce  its slight generalization: the category
of $i$--affine strict polynomial functors ${\cal P}_d^{af_i}$,hence  we start with 
reviewing its basic properties. Since  all the proofs from [C4] still work
in the present context,  the reader is referred for them to [C4].
The only exception where we provide a full proof is Proposition 3.2 which
was merely mentioned in [C4].\newline
Let $A_i:=\ka[x_1, x_2,\ldots,x_i]/(x_1^p, x_2^p,\ldots x_i^p)$ for $|x_j|=2p^j$
and let $\Gamma^d{\cal V}_{A_i}$ stands for the following graded $\ka$--linear category.  The  objects of $\Gamma^d{\cal V}_{A_i}$ are finite
dimensional vector spaces, though we follow the convention taken in [C4, Section 2] and label them as $V\ot A_i$ where $V$ is a finite dimensional vector space.
The morphisms are given as
\[\mbox{Hom}_{\Gamma^d {\cal V}_{A_i}}(V\ot A_i,W\ot A_i):=\Gamma^d(\mbox{Hom}(V,W)\ot A_i).\] 
An $i$--affine strict polynomial functor of degree $d$ is a graded functor from $\Gamma^d{\cal V}_{A_i}$ to the category of ${\bf Z}$--graded bounded below
finite dimensional in each degree vector spaces (c.f. [C4, Section 2]). The $i$--affine strict polynomial functors of degree $d$ form the $\ka$--linear 
graded abelian category ${\cal P}_d^{af_i}$ with morphisms being the natural transformations.  For any finite dimensional vector space $U$ we have the 
representable   $i$--affine strict polynomial functors of degree $d$
$h^{U\ot A_i}$ given by the formula
\[V\ot A_i\mapsto\mbox{Hom}_{\Gamma^d {\cal V}_{A_i}}(V\ot A_i,W\ot A_i)=\Gamma^d(\mbox{Hom}(V,W)\ot A_i),\]
and by the Yoneda lemma [C4, Prop. 2.2] we have
\[\mbox{Hom}_{{\cal P}_d^{af}}(h^{U\ot A_i},F)\simeq F(U\ot A_i).\]
Similarly, we have the co--representable functor $c_{U\ot A_i}^*$ given by
\[V\ot A\mapsto \mbox{Hom}_{\Gamma^d {\cal V}_{A}}(V\ot A,U\ot A)^*\]
where $(-)^*$ stands for the graded \ka--linear dual. This time the Yoneda lemma gives
\[\mbox{Hom}_{{\cal P}_d^{af}}(F,c_{U\ot A}^*)\simeq F(U\ot A)^*.\] 
Analogosly to the non--affine case, ${\cal P}_d^{af_i}$ is equivalent to some module category. Namely, let as define the $i$--affine Schur algebra $S_{d,n}^{af_i}:=
\Gamma^d(\mbox{End}(\ka^n))\ot A_i)$. Then $F(\ka^n\ot A_i)$ is naturally
a graded $S_{d,n}^{af_i}$--module and we have [C4, Prop.~2.5]
\begin{prop}
	If $n\geq d$ then
	\[ev_n: {\cal P}_d^{af_i}\ra S_{d,n}^{af_i}\mbox{-mod}^{f+},\]
	where $S_{d,n}^{af}\mbox{-mod}^{f+}$ is the category of bounded below finite dimensional in each degree graded $S_{d,n}^{af}$--modules,
	is an equivalence of graded abelian categories.
\end{prop}
The forgetful functor $z:\Gamma^d{\cal V}_{A_i}\ra \Gamma^d{\cal V}$
induces an exact functor $z^*:{\cal P}_d\ra {\cal P}_d^{af_i}$ which has
right and left adjoints $t^*,h^*:{\cal P}_d^{af_i}\ra {\cal P}_d$ (consult
[C4, Sect. 2] on grading issue). \newline Much deeper is relation between ${\cal P}_d^{af_i}$ and ${\cal P}_{dp^i}$, since it only emerges at the level of derived
categories.
In order to develop   homological algebra in  ${\cal P}_d^{af_i}$ we regard
 it as a DG category (with the trivial differential) (see [K1], [K2], [C4, Section 3]). Then we consider the category of complexes
over ${\cal P}_d^{af_i}$, i.e. the category of graded functors 
 from $\Gamma^d{\cal V}_{A_i}$ to the category of  bounded below
 complexes of finite dimensional in each degree vector spaces. The derived category ${\cal DP}_d^{af_i}$ is obtained from the category of complexes by
 inverting the class of quasiisomorphisms. This procedure can be conducted within the formalism of Quillen model categories. Namely the category
 of complexes over ${\cal P}_d^{af_i}$ can be equipped with either of two model structures: the projective
 one in which every object is fibrant and $h^{U\ot A_i}$ are cofibrant and the injective one in which every object is cofibrant and $c_{U\ot A_i}^*$ are fibrant.
 In both cases   the homotopy category is equivalent to  ${\cal DP}_d^{af_i}$.\newline
 The main result of [C4] is a construction of a full embedding
 \[{\bf C}^{af_i}:{\cal P}_d^{af_i}\ra{\cal P}_{dp^i}\]
 and its right adjoint 
 \[{\bf K}^{af_i}:{\cal P}_d^{af_i}\ra{\cal P}_{dp^i}\]
 called the affine derived Kan extension. This adjunction is compatible with the
 adjunction $\{{\bf C}, {\bf K^r}\}$  in the sense that we have isomorphisms
 of functors [C4, Theorem 5.1]:
 \[{\bf K}^r\simeq t^*\circ {\bf K}^{af},\ \ \ {\bf C}\simeq {\bf C}^{af}\circ z^*.\] 
 We finish our review by discussing the compatibility of  $\{{\bf C}^{af}, {\bf K^{af}}\}$ with the Kuhn duality. This is a non--trivial problem, which was  mentioned in [C4, Sect. 6]. In particular we shall use in the proof the main 
 result of [C3]. We need it in the present article  in order to connect the results of Sections 2
 and 5.\newline
 We recall that ${\cal DP}_d$ denotes the bounded derived category. We denote
 by  ${\cal DP}_d^+$ the derived category coming from the bounded below complexes.
 Let ${\cal DP}_d^{af_i,b}$  be the smallest full triangulated subcategory of ${\cal DP}_d^{af_i,}$ containing $h^{U\ot A_i}$ and closed under taking direct factors.
 In other words: ${\cal DP}_d^{af_i,b}$ is the full subcategory of ${\cal DP}_d^{af_i}$ consisting of all compact objects. 
 \begin{prop}
 We have the following isomorphisms of functors:
 \begin{enumerate}
 	\item \[(-)^{\#}\circ{\bf K}^{af}\circ (-)^{\#}\simeq {\bf K}^{af}\]
 	as functors between ${\cal DP}_{dp^i}^+$ and ${\cal DP}_{d}^{af}$,
 	\item \[(-)^{\#}\circ{\bf C}^{af}\circ (-)^{\#}\simeq {\bf C}^{af}\]
 	as functors between ${\cal DP}_d^{af_i,b}$ and ${\cal DP}_{dp^i}$.
 \end{enumerate}
 \end{prop}
 {\bf Proof:}  In fact the first part can be deduced from the proof of [C3, Theorem~2.1].
 Let $k: \Gamma^d{\cal V}_{A_i}\ra \Gamma^d{\cal V}$ be the functor induced by the projrction $A_i\ra \ka$ and let $k^*: {\cal DP}_d\ra {\cal DP}_{d}^{af_i,b}$
 be the functor induced by the precomposition with $k$. Then 
 \[{\bf K}^{af_i}(S^{dp^i}_{U^*})=k^*(S^d_{U^{(i)*}})\] and
 \[k^*(S^d_{U^{(i)*}})=k^*(\Gamma^d_{U^{(i)*}})[-2d(p^i-1)]\]
 c.f. [C4, Prop.~2.4.2]. Thus the isomorphism constructed in the proof of [C3, Theorem~2.1] can be interpreted as 
 \[{\bf K}^{af_i}(F^{\#})\simeq {\bf K}^{af_i}(F)^{\#}\]
 for any $F\in{\cal DP}_{dp^i}$. In order to extend this isomorphism to ${\cal DP}_{dp^i}^+$ we observe that, since ${\bf K}^{af_i}$ commutes with infinite direct colimits (because $P^{\bullet}$ is a compact object in ${\cal DP}_{dp^i}^+$), both left and right hand sides take direct colimits into codirect limits.\newline
 In order to obtain the second part we recall that 
 \[(h^{U\ot A_i})^{\#}=c^*_{U^*\ot A_i}[-2d(p^i-1)]=z^*(S^d_{U^*}).\]
 Hence
 \[{\bf C}^{af_i}((h^{U\ot A_i})^{\#})={\bf C}^{af_i}(z^*(S^d_{U^*}))={\bf C}(S^d_{U^*})=S^{d(i)}_{U^*}=\]
 \[(\Gamma^{d(i)}_{U^*})^{\#}=({\bf C}^{af_i}(h^{U\ot A_i}))^{\#},\]
 which gives the required isomorphism for any $F\in {\cal DP}_d^{af_i,b}$. This
 time, since ${\bf C}^{af}$ does not commute with infinite codirect limits, the isomorphim cannot be extended to the whole ${\cal DP}_d^{af_i}$. Indeed, as it was observed in [C4, Section~6], e.g. for $F=k^*(I)\in{\cal DP}_1^{af_1}$,  computing the both sides of the postulated isomorphism gives different results.\qed
 \subsection{Serre functor in ${\cal DP}_d^{af_i}$}
 In this subsection we introduce (a suitably modified version of) Serre functor
 in ${\cal DP}_d^{af_i}$. In fact a genuine Serre functor exisists only on 
 ${\cal DP}_d^{af_i,b}$ but this is not very useful for us since the objects
 of ${\cal DP}_d^{af_i}$ rarely have cofibrant replacements in 
 ${\cal DP}_d^{af_i,b}$. On the other hand one cannot hope for existing of a Serre
 functor in ${\cal DP}_d^{af_i}$, since it is not even a Hom--finite category.
 What we really have is the following weaker version of Serre functor:
 \begin{defi}
 	Let ${\cal C}$ be a \ka--linear category. A \ka--linear functor ${\bf S}:
 	{\cal C}\ra {\cal C}$ is a weak (left) Serre functor if:
 	\begin{enumerate}
 		\item There is a natural in $X,Y$ isomorphism
 		\[\mbox{Hom}_{{\cal C}}({\bf S}(X),Y)\simeq \mbox{Hom}_{{\cal C}}(Y,X)^*\]
 		whenever $X$ or $Y$ is compact.
 		\item ${\bf S}$ is an auto--equivalence.
 	\end{enumerate}
 \end{defi}
 For example, if ${\cal C}$ is the bounded derived category of category of finitely
 generated modules over a finite dimensional algebra of finite homological 
 dimension then by [BV] ${\cal C}$ posses a Serre functor. In that case, it is easy 
 to see that it extends to a weak Serre functor. However, our situation is quite different, since ${\cal P}_d^{af_i}$ is of infinite homological dimension.   
 The crucial property of
 ${\cal DP}_d^{af_i}$ which makes a construction analogous to that used in Section 2 working is the following technical fact:
 \begin{prop} For any finite dimensional space $U$,
 $c_{U\ot A_i}^*$ is a compact object of ${\cal DP}_d^{af_i}$, i.e. the functor
 $\mbox{Hom}_{{\cal DP}_d^{af_i}}(c_{U\ot A_i}^*,-)$ commutes with infinite 
 direct colimits.
 \end{prop} 
 Let $P^{\bullet}$ be a finite projective resolution of $S^d_{U^*}$ in ${\cal P}_d$.
 Then 
 \[h^*(P^{\bullet})\simeq h^*(S^d_{U^*})=c^*_{U\ot A_i}.\]
 Since $h^*\simeq t^*[2d(p^i-1)]$ by [C4, Prop.~2.4.5] and $t^*$ preserves cofibrant objects,
 $h^*(P^{\bullet})$ is a cofibrant replacement of   $c^*_{U\ot A_i}$. To conclude the proof we observe that since $P^{\bullet}$ is finite, $h^*(P^{\bullet})$ is compact. \qed
 Now we can define the affine Serre functor in a manner analogous to that used in Section  2. 
 \begin{defi}
 	We define a  functor $\ser: \dep_d^{af_i}\ra\dep_d^{af_i}$ by the formula
 	\[
 	{\bf S}^{af_i}(F)(V\ot A_i):=\mbox{Hom}_{{\cal DP}_d^{af_i}}(c_{V\ot A_i}^*,F).
 	\]
 \end{defi}
 We collect the basic properties of ${\bf S}^{af_i}$ which will be needed for the applications described in  Section 5.
\begin{theo}
The functor ${\bf S}^{af_i}$ satisfies the following properties:
\begin{enumerate}
\item There is a natural in $U\ot A_i$ isomorphism in 	${\cal P}_d^{af_i}$
\[{\bf S}^{af_i}(c^*_{U\ot A_i})\simeq h^{U\ot A_i}.\]
	\item  ${\bf S}^{af_i}$ is an auto--equivalence of ${\cal DP}_d^{af_i}$.
	\item ${\bf S}^{af_i}$ restricted to ${\cal DP}_d^{af_i,b}$ is a left Serre functor
	and it is a weak left Serre functor on the whole ${\cal DP}_d^{af_i}$.
	\item There are isomorphisms of functors
	\[{\bf S}^{af_i}\circ z^*[2d(p^i-1)]\simeq z^*\circ{\bf S},\ \ \ 
	{\bf S}\circ t^*\simeq t^*\circ{\bf S}^{af_i}[2d(p^i-1)],\]
	\[{\bf S}^{af_i}\circ {\bf K}^{af}\simeq {\bf K}^{af}\circ{\bf S},\ \ \
	{\bf S}\circ {\bf C}^{af_i}\simeq {\bf C}^{af}\circ{\bf S}^{af_i}.\] 
\end{enumerate}
\end{theo}
{\bf Proof:}  Since $c_{U\ot A_i}^*$ is fibrant in the injective Quillen structure
we get
\[{\bf S}^{af}(c_{U\ot A_i}^*)(V\ot A_i)\simeq
\mbox{Hom}_{{\cal P}_d^{af_i}}(c_{V\ot A_i}^*,c_{U\ot A_i}^*)\simeq
(c_{V\ot A_i}^*(U\ot A_i))^*=\]
\[\Gamma^d(\mbox{Hom}(U,V)\ot A_i)=h^{U\ot A_i}
(V\ot A_i).\]
In order to get the second part we define ``the right affine Serre functor'' ${\bf S}_r^{af_i}$ by the formula
\[
{\bf S}^{af_i}_r(F)(V\ot A_i):=\mbox{Hom}_{{\cal P}_d^{af_i}}(F,h^{V\ot A_i})^*.
\]
Then by the computation analogous to that giving the first part we show
that ${\bf S}_r^{af_i}(h^{U\ot A_i})=c_{U\ot A_i}^*$. 
Thus the transformation $id\ra {\bf S}^{af_i}\circ{\bf S}_r^{af_i}$ coming
from the Yoneda lemma is an isomorphsim for $h^{U\ot A_i}$. Hence 
${\bf S}^{af_i}$and ${\bf S}_l^{af_i}$ are mutually inverse on   ${\cal DP}_d^{af_i,b}$. Since, by Proposition 3.4, 
${\bf S}^{af_i}$ commutes with infinite direct colimits, it is an equivalence
on the whole ${\cal DP}_d^{af_i}$. \newline
In order to get the third part, again by Proposition 3.4, it suffices to establish a natural 
in $F\in{\cal DP}_d^{af_i}$ and $U\ot A_i$ isomorphism
\[\mbox{Hom}_{{\cal DP}_d^{af_i}}(F,c^*_{U\ot A_i})\simeq\mbox{Hom}_{{\cal DP}_d^{af_i}}({\bf S}^{af_i}(c^*_{U\ot A_i}),F)^*,\]
which by the first part and the fact that $h^{U\ot A_i}$ is cofibrant and $c^*_{U\ot A_i}$ is fibrant reduces to the isomorphism
\[\mbox{Hom}_{{\cal P}_d^{af_i}}(F,c^*_{U\ot A_i})\simeq\mbox{Hom}_{{\cal P}_d^{af_i}}(h^{U\ot A_i},F)^*,\]
which follows from the Yoneda lemma.\newline
In order to obtain the first isomorphism in part 4, we recall that
$z^*(\Gamma^d_{U^*})=h^{U\ot A_i}$ and $z^*(S^d_{U^*})=c^*_{U\ot A_i}[-2d(p^i-1)]$. Hence we get natural in $U$ isomorphisms:
\[{\bf S}^{af_i}\circ z^*(S^d_{U^*})={\bf S}^{af_i}(c^*_{U\ot A_i}[-2d(p^i-1)])=
h^{U\ot A_i}[-2d(p^i-1)]\]
and
\[z^*\circ{\bf S}(S^d_{U^*})=z^*(\Gamma^d_{U^*})=h^{U\ot A_i}.\]
The second isomorphism follows from the facts that
$t^*(c^*_{U\ot A_i})=S^d_{(U\ot A_i)^*}$ and $t^*(h^{U\ot A_i})=\Gamma^d_{(U\ot A_i)^*}[-2d(p^i-1)]$.\newline
The proof of the last two isomorphisms is analogous to that of Proposition 3.2.
The last formula holds on the whole ${\cal DP}_d^{af_i}$ because ${\bf C}^{af_i}$
commutes with infinite direct limits. \qed
{\bf Remark:} When we compose the isomorphisms from part 4, we obtain
the formulae:
\[{\bf S}\circ {\bf K}^{r}\simeq {\bf K}^{r}\circ{\bf S}[2d(p^i-1)],\ \ \
{\bf S}\circ {\bf C}\simeq {\bf C}\circ{\bf S}[-2d(p^i-1)],\] 
from which, in particular,  Theorem 2.2.3.1 follows. Thus we see that this shift phenomenon which produces the Poincar\'e duality  is  related to the scalar extension from ${\cal P}_d$
to ${\cal P}_d^{af_i}$. 
\section{Affine semiblocks}
In this section we introduce certain subcategories of ${\cal P}_d^{af_i}$ which
correspond to the blocks in ${\cal P}_d$. We call them semiblocks since they generate  ${\cal P}_d^{af_i}$ and we conjecture that ${\cal P}_d^{af_i}$ is stratified
by these subcategories. This structure may  be interesting for its own but in 
the present article  we are mainly interested in the Serre functor restricted to 
the semiblocks, since, as it will be shown in the next section, in  certain cases
it enjoys very special properties.\newline
We recall that the category ${\cal P}_d$ admits decmopsition into the blocks:
\[{\cal P}_d\simeq {\cal P}_{\la^1}\times\ldots\times {\cal P}_{\la^s}
\]
and the set of blocks is indexed by the family $\la^1,\ldots,\la^s$ of 
$p$--core Young diagrams of weight $d-pj$ for some $j\geq 0$. By the Yoneda
lemma, there is the corresponding decomposition of the bifunctor 
$(V,W)\mapsto\Gamma^d(\mbox{Hom}(W,V))$ into the ``block bifunctors'':
\[\Gamma^d(\mbox{Hom}(W,V))\simeq B_{\la^1}(V,W)\oplus\ldots B_{\la^s}(V,W).
\] 
The Cauchy decomposition [ABW, Th.~III.1.4] provides the filtration of bifunctor
$\Gamma^d(\mbox{Hom}(W,V))$  with the associated object
\[ \bigoplus_{\mu\in Y_d} W_{\mu}(V)\ot W_{\mu}(W^*)\]
where $Y_d$ stands for the set of Young diagrams of weight $d$.
Hence each $B_{\la}(V,W)$ has the filtration with the associated object
\[ \bigoplus_{\mu\in Y_{\la^j}} W_{\mu}(V)\ot W_{\mu}(W^*)\]
where $Y_{\la}$ is the set of Young diagrams of degree $d$ belonging to the
block labeled by $\la^j$.\newline
Moreover, the bifunctor $B_{\la^j}$ can be used to form the category
$B_{\la^j}{\cal V}$ whose objects are finite vector spaces and
\[
\mbox{Hom}_{B_{\la^j} {\cal V}}(V,W):=B_{\la^j}(W,V).
\]
Then the category ${\cal P}_{\la^j}$ can be identified with the category 
of $\ka$--linear functors from $B_{\la^j}{\cal V}$ to the category
of finite dimensional vector spaces over \ka.
The main objective of the present section is to define the affine counterpart
of ${\cal P}_{\la}$, relate it to ${\cal P}_{\la}$ and ${\cal P}_{dp^i}$,  and 
 equip it with a Serre functor.\newline
Let us fix a $p$--core Young diagram $\la$ of wieght $|\la|=d-pj$ and
let $B_{\la}^{(i)}$ denote the bifunctor $(V,W)\mapsto B_{\la}(V^{(i)},W)$. 
We introduce the graded category $B_{\la}{\cal V}_{A_i}$ with the objects being 
finite dimensional vector spaces and the morphisms given by the formula
\[\mbox{Hom}_{B_{\la}{\cal V}}(V\ot A_i,V'\ot A_i):=\mbox{Ext}^*_{{\cal P}_{dp^i}}(B_{\la}^{(i)}(-,V'),B_{\la}^{(i)}(-,V))\]
where we choose to label the objects by $V\ot A_i$
in order to make our terminology coherent with that used in Section 3 and [C4].
Thanks to the Collapsing Conjecture [C3, Cor. 3.7] the Hom--spaces in $B_{\la}{\cal V}_{A_i}$
admit a  more  explicit description. Namely, we have natural
in $V, V'$ isomorphisms
\[\mbox{Ext}^*_{{\cal P}_{dp^i}}(B_{\la}^{(i)}(-,V'),B_{\la}^{(i)}(-,V))\simeq
\mbox{Ext}^*_{{\cal P}_{d}}(B_{\la}(-,V'),B_{\la}(-\ot A_i,V))\simeq\]
\[B_{\la}(V'\ot A_i,V).\]
With this description the composition in $B_{\la}{\cal V}_{A_i}$ is given as the composite of the scalar extension:
\[
B_{\la}(V''\ot A_i,V')\ot B_{\la}(V'\ot A_i,V)\ra
B_{\la}(V''\ot A_i\ot A_i,V'\ot A_i)\ot B_{\la}(V'\ot A_i,V),
\]
the composition in $B_{\la}{\cal V}$:
\[ B_{\la}(V''\ot A_i\ot A_i,V'\ot A_i)\ot B_{\la}(V'\ot A_i,V)\ra B_{\la}(V''\ot A_i\ot A_i,V)\]
and the morphism 
\[
B_{\la}(V''\ot A_i\ot A_i,V)\ra B_{\la}(V''\ot A_i,V)\]
induced by the multiplication $A_i\ot A_i\ra A_i$. 
We then define ${\cal P}_{\la}^{af_i}$ as the category of graded $\ka$--linear functors from  $B_{\la}{\cal V}_{A_i}$ to the category 
of ${\bf Z}$--graded bounded below
finite dimensional in each degree vector spaces.\newline
The category ${\cal P}_{\la}^{af_i}$ shares with ${\cal P}_{d}^{af_i}$ its basic 
properties. In particular we have
 the representable functor $h_{\la}^{U\ot A_i}$ in ${\cal DP}_{\la}^{af_i}$ given explicitly by the formula
\[h_{\la}^{U\ot A_i}(V):=\mbox{Hom}_{B_{\la}{\cal V}_{A_i}}(U\ot A_i,V\ot A_i),
\]
 the corepresentable functor $c_{\la,U\ot A_i}^*$ and the block affine Kuhn duality. Also the analog 
of Proposition 3.1 holds.
Let us call the block
affine Schur alegbra the graded algebra 
\[\mbox{Hom}_{B_{\la}{\cal V}_{A_i}}(\ka^d\ot A_i,\ka^d\ot A_i)\simeq 
B_{\la}(\ka^d\ot A_i,\ka^d).\]
Then 
\begin{prop}
	The evaluation functor $F\mapsto F(\ka^d\ot A_i)$ gives an equivalence of graded categories
	\[{\cal P}_{\la}^{af_i}\simeq B_{\la}(\ka^d\ot A_i,\ka^d)-grmod.\]
\end{prop}
At last, the adjunction $\{z^*,t^*\}$ between ${\cal P}_d$ and ${\cal P}_d^{af_i}$
clearly extends to the adjunction $\{z^*_{\la},t^*_{\la}\}$ between ${\cal P}_{\la}$ and ${\cal P}_{\la}^{af_i}$.\newline
On the other hand, when we try to decompose  ${\cal P}_{d}^{af_i}$, into the  product of  ${\cal P}_{\la^j}^{af_i}$ we face a problem that for
$\la\neq \la'$, ${\cal P}_{\la}^{af_i}$ and ${\cal P}_{\la'}^{af_i}$ are not orthogonal
as subcategoeries of ${\cal P}_{d}^{af_i}$. We will come back to this observation 
later, since it is best comprehensible at the level of derived categories.\newline Now
we turn to describing relation between  ${\cal P}_{\la}^{af_i}$ and 
${\cal P}_{d}^{af_i}$ more precisely. 
Let \[i_{\la}: B_{\la}(V\ot A_i,W)\ra \Gamma^d(\mbox{Hom}(W,V\ot A_i))\]
be the natural embedding and 
\[\pi_{\la}: \Gamma^d(\mbox{Hom}(W,V\ot A_i))\ra B_{\la}(V\ot A_i,W)\]
be the natural projection. Then the composite $\epsilon_{\la}:=i_{\la}\circ \pi_{\la}$ can be thought of as an idempotent endofunctor on 
$\Gamma^d{\cal V}_{A_i}$ (being the identity on the objects). 
Thus the category $B_{\la}{\cal V}_{A_i}$ can be identified with the category
$\epsilon_{\la}(\Gamma^d{\cal V}_{A_i})\epsilon_{\la}$ whose objects are those
of $\Gamma^d{\cal V}_{A_i}$ but
\[\mbox{Hom}_{\epsilon_{\la}(\Gamma^d{\cal V}_{A_i})\epsilon_{\la}}(V,V'):=
\epsilon_{\la}(\mbox{Hom}_{\Gamma^d{\cal V}_{A_i}}(V,V'))\epsilon_{\la}.\]
Then the assignement \[(V,V')\mapsto \epsilon_{\la}(\mbox{Hom}_{\Gamma^d{\cal V}_{A_i}}(V,V'))\]
defines a $\Gamma^d{\cal V}_{A_i}\times B_{\la}{\cal V}_{A_i}$--bimodule
in the terminoloy of [K1, Sect. 6]. Hence we get a
 pair of functors $j^*_{\la}, j_{\la*}$ which satisfy the following properties. 
\begin{prop}
	\begin{enumerate}
		\item	The functor $j_{\la*}$ is right adjoint to $j^*_{\la}$.
		\item  The functor 
		$j_{\la}^*:  {\cal P}_{\la}^{af_i}\ra {\cal P}_d^{af_i}$
		is a full embedding.
	\end{enumerate}
\end{prop}
{\bf Proof: } The adjunction follows from the machinery of standard functors
developed in [K1, Sect. 6]. The full embedding follows from the fact that 
$j^*_{\la}\circ j_{\la*}\simeq id$. \qed
Let us remark that Proposition 4.2 may be thought of as a categorification of
[CPS, Prop. 2.1]. This explains our choice of notations with $j^*, j_*$ instead of $H_X, T_X$ used in [K1].
In fact we could derive Proposition 4.2 directly from [CPS, Prop. 2.1] by invoking
our Proposition 4.1 but we prefer to  consistently work in functor categories.
We also mention that Proposition 4.2 carries over to the level of derived categories
which was the main objective of [K1] and [CPS] and which will be discussed in the next paragraph.\newline
Namely, 
we define     ${\cal DP}_{\la}^{af}$ as the derived category of DG--category
${\cal P}_{\la}^{af_i}$ in the manner analogous to that in Section 3. The
adjunctions $\{z^*_{\la},t^*_{\la}\}$ and $\{j^*_{\la},j_{\la*}\}$ carry over to the derived categories and, as we have already mentioned, the analog of the second part 
of Proposition 4.2 holds, i.e. we have the full embedding 
\[j_{\la}^*:  {\cal DP}_{\la}^{af_i}\ra {\cal DP}_d^{af_i},\]
which allows us to regard ${\cal DP}_{\la}^{af_i}$ as a full subcategory
of ${\cal DP}_d^{af_i}$. Then it is clear that our construction is compatible with the scalar extension from ${\cal DP}_d$ to ${\cal DP}_d^{af_i}$:
\begin{prop}
	There are isomorphisms of functors:
	\[ t^*\circ j^*_{\la}\simeq b_{\la}^*\circ t^*,\ \ \ z^*\circ b_{\la*}\simeq 
	j_{\la*}\circ z^*,\]
where $b_{\la}^*$ and $b_{\la*}$ are induced repsectively by the embedding of and the projection onto  the block.
\end{prop}
 
Now we would like to to construct a block version of the affine derived Kan extension in order to relate ${\cal DP}_{\la}^{af_i}$ to ${\cal DP}_{p^id}$.
For this we need an 
 analog of the formality result [C4, Th.~4.2].
Let $X_{\la}$ be a projective
resolution of $B_{\la}^{(i)}$ in ${\cal P}_{dp^i}^d$. We introduce a DG category
$\Gamma^d {\cal V}_{X_{\la}}$ with the objects being finite dimensional
vector spaces and 
\[\mbox{Hom}_{\Gamma^d{\cal V}_{X_{\la}}}(V,V'):=\mbox{Hom}_{{\cal P}_{dp}}(X_{\la}(-,V'),X_{\la}(-,V)).\]
Then $B_{\la}{\cal V}_{A_i}$ is clearly the cohomology category of $\Gamma^d {\cal V}_{X_{\la}}$ but we have a much stronger result (c.f. [C4, Th.~4.2]):
\begin{prop}
The identity on the objects extends to an equivalence of  DG categories
$\phi_{\la}: B_{\la}{\cal V}_{A_i}\simeq \Gamma^d {\cal V}_{X_{\la}}$.
\end{prop}
{\bf Proof:}  
Since $\Gamma^d(I^{(i)}\ot I^*)\simeq B^{(i)}_{\la}\oplus B'$, we can obtain $X$, the projective resolution of
$\Gamma^d(I^{(i)}\ot I^*)$, as the direct sum $X=X_{\la}\oplus X'$ of projective
resolutions of $B^{(i)}_{\la}$ and  $B'$.  Let 
\[i_{\la}: \mbox{Ext}^*_{{\cal P}_{dp}}(B_{\la}^{(i)}(-,V'),B_{\la}^{(i)}(-,V))\ra
\mbox{Ext}^*_{{\cal P}_{dp}}(\Gamma^d((-)^{(i)},V'),\Gamma^d((-)^{(i)},V))\]
be the embedding induced by the decomposition $\Gamma^d(I^{(i)}\ot I^*)\simeq B^{(i)}_{\la}\oplus B'$ 
(we have already encountered this embedding when constructing the idempotent functor $\epsilon_{\la})$. Let us define similarly the projection
\[\widetilde{\pi}_{\la}: \mbox{Hom}_{{\cal P}_{dp}}(X(-,V'),X(-,V))\ra
\mbox{Hom}_{{\cal P}_{dp}}(X_{\la}(-,V')X_{\la}(-,V)).\]
We define $\phi_{\la}: B_{\la}{\cal V}_{A_i}\ra\Gamma^d {\cal V}_{X_{\la}}$
as the composite
$\phi_{\la}:=\widetilde{\pi}_{\la}\circ\phi\circ i_{\la}$ where $\phi: \Gamma^d{\cal V}_{A_i}
\ra \Gamma^d{\cal V}_X$ is the transformation from [C4, Theroem 4.2] or rather
its multitwist analog (in fact this generalization is not entirely trivial, we refer the reader to [C5, Theorem~3.1] where an analogous construction is conducted in even greater generality). Then  the fact that $\phi_{\la}$ is an
equivalence easily follows from the fact that $\phi$ is an equivalence and that it
commutes with the idempotent  $\epsilon_{\la}:=i_{\la}\circ \pi_{\la}$
and its $\Gamma^d {\cal V}_{X_{\la}}$--analog $\widetilde{\epsilon}_{\la}:=\widetilde{i}_{\la}\circ \widetilde{\pi}_{\la}$. \qed
Thanks to Proposition 4.4 we are able to construct the block affine derived Kan extension.
We summarize its basic properties below
\begin{prop}
	There exist functors ${\bf C}^{af_i}_{\la}: {\cal DP}_{\la}^{af_i}\ra {\cal DP}_{dp^i}$
	and ${\bf K}^{af_i}_{\la}: {\cal DP}_{dp^i}\ra {\cal DP}_{\la}^{af_i}$ satisfying the
	following properties:
	\begin{enumerate}
		\item ${\bf K}^{af_i}_{\la}$ is right adjoint to ${\bf C}^{af_i}_{\la}$.
		\item ${\bf C}^{af_i}_{\la}$ is a full embedding.
		\item The functors ${\bf C}^{af_i}_{\la}$ (restricted to
		${\cal DP}_{\la}^{af_i,b}$) and ${\bf K}^{af_i}_{\la}$ commute with the Kuhn duality. 
		\item There are isomorphisms of functors:
		\[ {\bf C}^{af_i}\circ j_{\la*}\simeq {\bf C}^{af_i}_{\la},\ \ \ 
		j_{\la}^*\circ {\bf K}^{af_i}\simeq {\bf K}^{af_i}_{\la}.
		\]
	\end{enumerate}
\end{prop}
The proofs of parts 1, 2, 3 are analogous to those of [C4, Th.~5.1] and our Proposition 3.2. The compatibilty formula in part 3 follows immediately from the constuction of the considered functors. \qed
Having at our disposal the block affine drived Kan extension we can offer 
a better explanation of the phenomenon of non--orthogonality of semiblocks. Namely let us 
take $F\in{\cal P}_{\la}$, $G\in{\cal P}_{\la'}$ for $\la\neq \la'$. Then
\[\mbox{Hom}^*_{{\cal DP}_d^{af_i}}(j^*_{\la}(z^*_{\la}(F)),j^*_{\la'}(z^*_{\la'}(G)))\simeq
\mbox{Hom}^*_{{\cal DP}_{dp^i}}({\bf C}^{af_i}(j^*_{\la}(z^*_{\la}(F))), {\bf C}^{af_i}(j^*_{\la'}(z^*_{\la'}(G))))\simeq\] 
\[\simeq\mbox{Ext}^*_{{\cal P}_{dp^i}}(F^{(i)},G^{(i)})\]
and the latter Ext--groups may well be non--trivial. 
Thus we see that the reason for the non--orthogonality of semiblocks is simply that 
 the Frobenius twist transfers all the blocks from ${\cal P}_d$
into the single (pricipal) block in ${\cal P}_{dp^i}$. Still, the block decomposition on ${\cal P}_d$ when pushed to  ${\cal P}_d^{af_i}$ generates certain structure
on the latter category. It seems that we have on ${\cal P}_d^{af_i}$ (and
its derived category) a stratification with the set of strata indexed by the blocks
of ${\cal P}_d$ but we defer a deeper study of semiblock structure to a future work.\newline  
We finish this section  by  endowing the category ${\cal DP}_{\la}^{af_i}$ with a Serre functor. This may be achieved by a construction analogous to that given in the global (affine) case.  
\begin{defi}
We define the block affine Serre functor ${\bf S}^{af_i}_{\la}:{\cal DP}_{\la}^{af_i}
\ra{\cal DP}_{\la}^{af_i}$ by the formula
\[
{\bf S}^{af_i}_{\la}(F)(V\ot A_i):=\mbox{RHom}_{{\cal P}_{\la}^{af_i}}(c_{\la,V}^*,F).
\]
\end{defi}
Then the block analog of Proposition 3.5 holds
\begin{theo}
The functor ${\bf S}^{af_i}_{\la}$ satisfies the following properties:
\begin{enumerate}
	\item There is a natural in $U\ot A_i$ isomorphism in 	${\cal P}_{\la}^{af_i}$
	\[{\bf S}^{af_i}_{\la}(c^*_{\la,U\ot A_i})\simeq h^{U\ot A_i}_{\la}.\]
	\item  ${\bf S}^{af_i}_{\la}$ is an auto--equivalence of ${\cal DP}_d^{af_i}$.
	\item ${\bf S}^{af_i}_{\la}$ restricted to ${\cal DP}_{\la}^{af_i,b}$ is a left Serre functor
	and it is a weak left Serre functor on the whole ${\cal DP}_{\la}^{af_i}$.
	\item There are isomorphisms of functors
	\[{\bf S}^{af_i}_{\la}\circ z^*_{\la}[2d(p^i-1)]\simeq z^*_{\la}\circ{\bf S}_{\la},
	\ \ \ 
	{\bf S}_{\la}\circ t^*_{\la}\simeq t^*_{\la}\circ{\bf S}^{af_i}_{\la}[2d(p^i-1)],\]
	\[{\bf S}^{af_i}_{\la}\circ {\bf K}^{af}_{\la}\simeq {\bf K}^{af}_{\la}\circ
	{\bf S}_{\la},\ \ \
	{\bf S}_{\la}\circ {\bf C}^{af_i}_{\la}\simeq {\bf C}^{af}_{\la}\circ
	{\bf S}^{af_i}_{\la},\] 
where ${\bf S}_{\la}$ is ${\bf S}$ restricted to the block ${\cal DP}_{\la}$.
\end{enumerate}
\end{theo}

The proof of Theorem 3.6 carries over to the current situation.\section{Calabi--Yau structure on basic affine blocks}
In this section we show that the affine Serre functor when restricted to certain
semiblocks in ${\cal DP}_d^{af_i}$ is isomorphic to the shift functor.\newline 
We recall that  a block in ${\cal P}_d$ is called basic if it contains a single simple 
object. Hence the basic blocks are indexed by $p$--core Young diagrams of weight
$d$ and we also call such Young diagrams basic. So, let us fix  a basic Young diagram $\lambda$. 
Then $S_{\la}\simeq W_{\la}\simeq F_{\la}$. Moreover, $S_{\la}$ is injective and projective
and every object of $P_{\la}$ is a direct sum of $S_{\la}$, therefore the category ${\cal P}_{\la}$ is semisimple.\newline
We recall that a triangulated category ${\cal T}$ with a Serre functor ${\bf S}_{{\cal T}}$ is called {\em Calabi--Yau of dimension $n$} if there is an isomorphism of
functors ${\bf S}_{{\cal T}}\simeq \mbox{id}[n]$. Then we call a triangulated category ${\cal T}$ {\em weak Calabi--Yau of dimension $n$}  if it has a weak Serre functor ${\bf S}_{\cal T}$
such that ${\bf S}_{{\cal T}}\simeq \mbox{id}[n]$.
\begin{theo}
For any  basic Young diagram $\la$, the category ${\cal DP}_{\la}^{af_i,b}$
is Calabi--Yau of dimension $2d(p^i-1)$, the category
${\cal DP}_{\la}^{af_i}$
is weak Calabi--Yau of dimension $2d(p^i-1)$.
\end{theo} 
{\bf Proof: } The theorem is a formal consequence of the following properties
of the bifunctor $B_{\la}$ (the crucial second property is specific to basic blocks).
\begin{lem}
	There are the following isomorphisms of bifunctors:
	\begin{enumerate}
\item $B_{\la}(V,W\ot A_i)\simeq B_{\la}(V\ot A_i^*,W)$ for any $\la$.
\item 	$B_{\la}(V,W)\simeq S_{\la}(V)\ot S_{\la}(W^*)$ for basic $\la$.
	\end{enumerate}
\end{lem}
{\bf Proof of the Lemma} We recall that 
\[B_{\la}(V,W)=\mbox{Hom}_{{\cal P}_d}(B_{\la}(-,W),B_{\la}(-,V))\]
and a general fact that
\[\mbox{Hom}_{{\cal P}_d}(F(-\ot X),G)\simeq \mbox{Hom}_{{\cal P}_d}(F,G(-\ot X^*))\]
for any graded space $X$ and $F,G\in{\cal P}_d$.  This gives the first isomorphism.\newline
The second isomorphism immediately follows from the Cauchy decomposition
and the fact that $S_{\la}\simeq W_{\la}$ for basic $\la$. \qed
We recall that we deal with left Serre functors, hence we should show that
${\bf S}^{af_i}_{\la}\simeq id[-2d(p^i-1)]$.\newline
Since 
\[
{\bf S}^{af_i}_{\la}(F)(V\ot A_i):=\mbox{RHom}_{{\cal P}_{\la}^{af_i}}(c_{\la,V}^*,F)
\]
and by the Yoneda lemma
\[
F(V\ot A_i):=\mbox{RHom}_{{\cal P}_{\la}^{af_i}}(h^V_{\la},F),
\]
it suffices to find a natural in $V$ isomorphism
\[c^*_{\la,V}\simeq h^V_{\la}[2d(p^i-1)].\]
On one hand we have:
\[h^V_{\la}(W)=\mbox{Hom}_{B_{\la}{\cal V}_{A_i}}(V,W)=(B_{\la}(W\ot A_i,V))\simeq
S_{\la}(W\ot A_i)\ot S_{\la}(V^*),\]
on the other hand: 
\[c^*_{\la,V}(W)=(\mbox{Hom}_{B_{\la}{\cal V}_{A_i}}(W,V))^*=
(B_{\la}(V\ot A_i,W))^*\simeq (B_{\la}(V,W\ot A_i^*))^*\simeq\]
\[(S_{\la}(V)\ot S_{\la}(W^*\ot A_i))^*=
S_{\la}(V^*)\ot S_{\la}(W\ot A_i^*).
\]
Since $A_i^*\simeq A_i[2(p^i-1)]$ we have an isomorphism of functors
\[S_{\la}(-\ot A_i^*)\simeq S_{\la}(-\ot A_i)[2d(p^i-1)]\]
which completes the proof. \qed
\begin{cor}
	The category ${\cal DP}_1^{af_i,b}$ is Calabi--Yau of dimension $2(p^i-1)$
	and ${\cal DP}_1^{af_i}$ is weak Calabi--Yau of dimension $2(p^i-1)$.
\end{cor}
{\bf Proof: } The corollary follows from Theorem 5.1 and the fact that ${\cal P}_1$ consists
of a single block which is obviously basic.\qed
This fact has the following global generalization.
\begin{prop}
	For any $d<p$, 	the category ${\cal DP}_d^{af_i,b}$ is Calabi--Yau of dimension $2d(p^i-1)$
	and ${\cal DP}_d^{af_i}$ is weak Calabi--Yau of dimension $2d(p^i-1)$.
\end{prop}
{\bf Proof:} In fact for $d<p$ all the blocks in ${\cal P}_d$ are basic but since 
${\cal P}_d^{af_i}$ is not a product of its affine semiblocks, our statement cannot
be directly deduced from Theorem 5.1.  Instead one can repeat the proof of Theorem 5.1 in the present context. The crucial fact is that $\Gamma^d\simeq S^d$ if $d<p$. We leave the straightforward details to the reader. \qed
As we have said in the Introduction, the Calabi--Yau structure on ${\cal DP}^{af_i}_{\la}$
provides sort of categorical interpretation of the Poincar\'e duality. 
Hence it is not surprising that  one can deduce Corollary 2.4 from
Theorem 5.1 (and the compatibilty of the (block ) affine derived Kan extension with 
the Kuhn duality).\newline  Namely, by the block affine derived Kan extension we obtain
\[\mbox{Ext}^s_{{\cal P}_{p^id}}(F_{\mu}^{(i)}, F_{\lambda})\simeq
\mbox{Hom}_{{\cal DP}_{\la}^{af_i}}(z_{\la}^*(F_{\la}), {\bf K}^{af_i}_{\la}(F_{\mu})[s]).\]
Then we apply the Calabi--Yau isomorphism (we emphasize that we need ``the weak Calabi--Yau structure'' here, since ${\bf K}^{af_i}$ does not preserve compact objects)
\[\mbox{Hom}_{{\cal DP}_{\la}^{af_i}}(z_{\la}^*(F_{\la}), {\bf K}^{af_i}_{\la}(F_{\mu})[s])\simeq
\mbox{Hom}_{{\cal DP}_{\la}^{af_i}}( {\bf K}^{af_i}_{\la}(F_{\mu})[s-2d(p^i-1)], z_{\la}^*(F_{\la}))^*.\]
Next we apply the Kuhn duality and use the fact that it commutes with $z^*$ and 
${\bf K}^{af_i}_{\la}$
\[
\mbox{Hom}_{{\cal DP}_{\la}^{af_i}}( {\bf K}^{af_i}_{\la}(F_{\mu})[s-2d(p^i-1)], z_{\la}^*(F_{\la}))^*\simeq 
\mbox{Hom}_{{\cal DP}_{\la}^{af_i}}(z_{\la}^*(F_{\la}^{\#})[s-2d(p^i-1)], {\bf K}^{af_i}_{\la}(F_{\mu}^{\#}))^*.
\]
At last we come back to ${\cal DP}_{dp^i}$:
\[
\mbox{Hom}_{{\cal DP}_{\la}^{af_i}}(z_{\la}^*(F_{\la}^{\#})[s-2d(p^i-1)], {\bf K}^{af_i}_{\la}(F_{\mu}^{\#}))^*\simeq
\mbox{Hom}_{{\cal DP}_{dp^i}}(F_{\la}^{(i)\#}[s-2d(p^i-1)], F_{\mu}^{\#}))^*\]
	and by using selfduality of simples we finally obtain our formula
\[
\mbox{Hom}_{{\cal DP}_{dp^i}}(F_{\la}^{(i)\#}[s-2d(p^i-1)], F_{\mu}^{\#}))^*\simeq	
\mbox{Ext}^{2d(p^i-1)-s}_{{\cal P}_{dp^i}}(F_{\la}^{(i)}, F_{\mu})^*.\]
Of course this approach is technically much more involved than that taken in Section 2, but  it shows how classical and affine phenomena are related and
also explains why we insist on considering weak Serre functors.\newline   
We finish our paper  by providing   description of  ${\cal P}_{\la}^{af_i}$ as
 a category of graded modules over certain explicitly described
graded algebra. Of course, for any $\la$, the category ${\cal P}_{\la}^{af_i}$ is equivalent to the category of graded modules over the block affine Schur algebra
by Proposition 4.1,
but this fact  is  not very useful in practice, since this  graded algebra is quite complicated. 
However, in the case of basic block the situation massively simplifies. First of all,
as we observed in Lemma 5.2 we have an isomorphism of graded vector spaces
\[B_{\la}(\ka^d\ot A_i,\ka^d)\simeq S_{\la}(\ka^d\ot A_i)\ot S_{\la}(\ka^{d*}).\]
However, in order to understand the multplicative structure it is better to take
a bit diiferent point of view.
Namely, by  Lemma 5.2  we have decomposition
\[
B_{\la}(-,\ka^d)\simeq\bigoplus_{j=1}^{s_{\la,d}}S_{\la}\] 
where $s_{\la,d}=\dim(S_{\la}(\ka^d))$.
Let us define a graded algebra
\[A_{i,\la}:=\mbox{Ext}^*_{{\cal P}_{dp^i}}(S_{\la}^{(i)},S_{\la}^{(i)}).\]
Then we have isomorphisms of graded algebras
\[B_{\la}(\ka^d\ot A_i,\ka^d)\simeq 
\mbox{Ext}^*_{{\cal P}_{dp^i}}(B^{(i)}_{\la}(-,\ka^d),B^{(i)}_{\la}(-,\ka^d))\simeq
\mbox{Ext}^*_{{\cal P}_{dp^i}}(\bigoplus_{j=1}^{s_{\la,d}}S_{\la}^{(i)},
\bigoplus_{j=1}^{s_{\la,d}}S_{\la}^{(i)})\simeq\]
\[M_{s_{\la,d}}(A_{i,\la}).\]
Since any matrix algebra is Morita equivalent to the ground algebra, we obtain
\begin{prop}
For any  basic Young diagram $\la$, the categories ${\cal P}_{\la}^{af_i}$
and $A_{i,\la}-grmod$ are equivalent.
\end{prop}
At last, let us take a look at the graded algebra $A_{i,\la}$. 
Firstly,  by the Collapsing Conjecture
\[A_{i,\la}\simeq \mbox{Hom}_{{\cal P}_d}(S_{\la}, S_{\la, A_i}).
\]
The dimension of the latter algebra can be explicitly expressed in terms of
the Littlewood--Richardson numbers. This point of view also allows one to 
describe the multiplication: it comes as the composite of scalar extension, Hom--multiplication
and multiplication in $A_i$:
 \[\mbox{Hom}_{{\cal P}_d}(S_{\la}, S_{\la, A_i})\ot \mbox{Hom}_{{\cal P}_d}(S_{\la}, S_{\la, A_i})\ra \]
 \[\mbox{Hom}_{{\cal P}_d}(S_{\la}, S_{\la, A_i})\ot \mbox{Hom}_{{\cal P}_d}
 (S_{\la\ot A_i}, S_{\la, A_i\ot A_i})\ra
 \mbox{Hom}_{{\cal P}_d}(S_{\la}, S_{\la, A_i\ot A_i})
 \ra\]
 \[\mbox{Hom}_{{\cal P}_d}(S_{\la}, S_{\la, A_i}).\]
 A bit different description of $A_{i,\la}$ is perhaps even more down to earth. It
follows from the fact that since $S_{\la}$ is a direct summand in $I^d$, there exists an idempotent $e_{\la}\in\ka[\Sigma_d]$ such that $S_{\la}=e_{\la} I^d$. Therefore we get
\[A_{i,\la}\simeq e_{\la} (\mbox{Ext}^*_{{\cal P}_{dp^i}}(I^{d(i)},I^{d(i)}))e_{\la}\simeq
e_{\la}(A_i^{\ot d}\ot\ka[\Sigma_d])e_{\la}.
\]
A subtle point here is that even if we would take the whole $S_{\la}$--isotypical summand in $I^d$ and the corresponding central idempotent $e_{\la}'$,
this $e_{\la}'$ is not central in the algebra $A_i^{\ot d}\ot\ka[\Sigma_d]$.
Hence $A_{i,\la}$ is not Morita equivalent to a direct factor in $A_i^{\ot d}\ot\ka[\Sigma_d]$.  This is another manifestation of the fact that 
affine semiblocks are not genuine blocks.\newline
All these descriptions drastically simplify for $d=1$. 
In this case we just obtain 
\begin{cor}
The categories ${\cal P}_{1}^{af_i}$
and $A_i-grmod$ are equivalent.
\end{cor}

\end{document}